\documentclass[a4paper,12pt]{amsart}

\usepackage{amssymb}
\usepackage{latexsym} 
\usepackage{amsfonts} 
\usepackage{amsmath}
\usepackage{eucal} 
\usepackage{bm} 
\usepackage{bbm} 
\usepackage{graphicx} 
\usepackage[english]{varioref} 
\usepackage[nice]{nicefrac} 
\usepackage[all]{xy}
\usepackage{amsthm}

\newcommand{\End}{{\rm End}}

\newcommand{\Spec}{\textrm{Spec}}

\newcommand{\Adm}{\textrm{Adm}}

\def\i{^{-1}}

\def\le{\leqslant}
\def\<{\langle} 
\def\>{\rangle}

\def\a{\alpha}

\def\g{\gamma}
\def\G{\Gamma}

\def\L{\Lambda}

\def\l{\lambda}

\def\ZZ{\mathbb Z}

\def\QQ{\mathbb Q}
\def\FF{\mathbb F}

\def\PP{\mathbb P}

\def\ca{\mathcal A}

\def\cg{\mathcal G}

\def\ci{\mathcal I}

\def\ck{\mathcal K}

\def\co{\mathcal O}
\def\cp{\mathcal P}

\def\ub{\underline B}
\def\up{\underline P}
\def\uh{\underline H}

\def\Fl{\mathcal Fl}

\def\tW{\tilde W}

\def\Spec{\rm Spec}
\def\loc{\rm loc}

\theoremstyle{plain}
\newtheorem{thm}{Theorem}[section] 
\newtheorem*{thm*}{Theorem} 
 \newtheorem{prop}[thm]{Proposition}
 
 \newtheorem{cor}[thm]{Corollary}

\theoremstyle{definition}

\theoremstyle{remark}

\newtheorem*{rmk}{Remark}
\newtheorem*{claim*}{Claim}

\begin{document}
\author{Xuhua He}
\address{Department of Mathematics, The Hong Kong University of Science and Technology, Clear Water Bay, Kowloon, Hong Kong}
\email{maxhhe@ust.hk}
\thanks{The author is partially supported by HKRGC grant 602011.}
\title[]{Normality and Cohen-Macaulayness of local models of Shimura varieties}
\keywords{Shimura variety, local model, affine flag, wonderful compactification}
\subjclass[2000]{14G35, 14M15, 14M27}

\begin{abstract}
We prove that in the unramified case, local models of Shimura varieties with Iwahori level structure are normal and Cohen Macaulay. 
\end{abstract}

\maketitle

\section*{Introduction}

Local models of Shimura varieties are projective schemes over the spectrum of a discrete valuation ring. Their singularities are expected to model the singularities that arise in the reduction modulo $p$ of Shimura varieties, with parahoric level structure. Local models also appear in the study of singularities of other moduli schemes (see Faltings \cite{F} and Kisin \cite{K}). We refer to the survey article by Pappas, Rapoport and Smithling \cite{PRS} for more details.

The simplest case of local models is for the modular curve with $\G_0(p)$-level structure. In this case, the local model is obtained by blowing up the projective line $\PP_{\ZZ_p}^1$ over $\Spec(\ZZ_p)$ at the origin of the special fiber $\PP^1_{\FF_p}=\PP^1_{\ZZ_p} \times_{\Spec \ZZ_p} \Spec \FF_p$. 

More generally, local models of Shimura varieties of PEL type with parahoric level structure were given by Rapoport and Zink in \cite{RZ} and in the ramified PEL case, by Pappas and Rapoport \cite{PR1}, \cite{PR2} and \cite{PR4}. The constructions there are representation-theoretic and mostly done case-by-case. 

Very recently, Zhu \cite{Z} (for equal characteristic analogy),  Pappas and Zhu \cite{PZ} made some new progress in the study of local models. They provide a group theoretic definition of local models that is not tied to a particular representation. The local model is constructed starting from the ``local Shimura data'' $(G, K, \{\mu\})$, where $G$ is a connected reductive group over $\QQ_p$, $K \subset G(\QQ_p)$ is a parahoric subgroup and $\{\mu\}$ is a geometric conjugacy class of one-parameter subgroups of $G$. Assume furthermore that $G$ splits over a tamely ramified extension of $\QQ_p$ and $\mu$ is minuscule. In \cite[Definition 7.1]{PZ}, Pappas and Zhu defined the local model $M^{\loc}$, which is a flat, projective scheme over $\Spec(\co_E)$. Here $E$ is the field of definition of $\mu$, $\co_E$ is the ring of integers of $E$ and $k_E$ its residue field. 

It is conjectured in \cite{PZ} that $M^{loc}$ is normal and Cohen-Macaulay. This question is also asked by Pappas, Rapoport and Smithling in \cite{PRS}. In this paper, we'll show that $M^{loc}$ is normal and Cohen-Macaulay in the unramified case. The precise statement will be found in Theorem \ref{main}. 

Now we discuss the outline of the proof. The generic fiber of $M^{\loc}$ is easy to understand. It is the Grassmannian variety associated to $\mu$. The special fiber $M^{\loc} \otimes_{\co_E} k_E$, on the other hand, is much more difficult to understand. A basic technique, introduced by G\"ortz \cite{G1}, is to embed the special fiber into an appropriate affine flag variety. 

One of the main results in \cite{PZ} is that the special fiber is the reduced union of affine Schubert varieties in the affine flag variety, indexed by the $\mu$-admissible set $\Adm^K(\mu)$ of Kottwitz and Rapoport, and each irreducible component of the special fiber is normal and Cohen-Macaulay. This extends results of G\"ortz \cite{G1}, \cite{G2}, Pappas and Rapoport \cite{PR1}, \cite{PR2}, \cite{PR4} on some Shimura varieties of PEL type. It is a deep result, based on the geometry of Schubert varieties in affine flag varieties \cite{F1} and \cite{PR3} and the coherence conjecture of \cite{PR3} recently proved by Zhu in \cite{Z}. 

By Serre's normality criterion and the behavior of depth under flat morphisms (see \cite[Remark 2.1.3]{PRS}), it remains to prove that the special fiber, as a whole, is Cohen-Macaulay. This is what we are going to do in this paper. 

The statement is obvious when the special fiber is irreducible (see e.g. \cite[Corollary 8.6]{PZ}). The main difficulty appears when the special fiber has more than one irreducible component. In \cite{G1}, G\"ortz proposes a combinatorial approach to this question and verifies the statement for unramified unitary groups of rank $\le 6$ in this way (with the aid of computer). Later, he checked a few more cases for $GL_n$ and $GSp_{2n}$ with small $n$ (unpublished). 

Our method here is quite different and more geometric. We now explain our strategy in more details. For simplicity, we only discuss the case where $G$ is split, of adjoint type and $K$ is an Iwahori subgroup of $G$. Let $\bar k$ be an algebraic closure of $\FF_p$. Let $LG=G_{\bar k((u))}$ be the loop group and $\ci$ be an Iwahori subgroup of $LG$. The projection map $G({\bar k[[u]]}) \to G_{\bar k}$ sends $\ci$ to a Borel subgroup $B$ of $G_{\bar k}$. Let $\Fl=LG/\ci$ be the affine flag variety. 

The main idea is to relate the affine flag variety with the wonderful compactification $X$ \cite{CP} of $G_{\bar k}$ and the geometric special fiber $M^{\loc} \otimes_{\co_E} \bar k$ with the boundary in $X$ of the parabolic subgroup of $G$ associated to $\mu$. 

The idea of relating the affine flag variety with the wonderful compactification comes from Springer. In \cite{S2}, Springer introduced a map from the loop group $LG$ to $X$, which factors through $G_{\bar k((u))}/\ck_1 \to X$. Here $\ck_1$ is the kernel of the projection map $G({\bar k[[u]]}) \to G_{\bar k}$. Notice that the natural map $LG/\ck_1 \to \Fl$ is a $B$-torsor. This map was used later by the author in \cite{H2} in the study of affine Deligne-Lusztig varieties in affine flag varieties. 

Springer's map is not continuous since the inverse image of the open subvariety $G_{\bar k}$ of $X$ is $G({\bar k[[u]]})/\ck_1$, which is not open in $LG/\ck_1$. However, as we'll see in Proposition \ref{11}, its restriction $$G({\bar k[[u]]}) s_{\l} G({\bar k[[u]]})/\ck_1 \to X$$ is a morphism. Here $\l$ is a coweight and $s_\l$ is the associated point in $G_{\bar k((u))}$. 

In particular, for the minuscule coweight $\mu$, we have the following diagram
\[\xymatrix{& \tilde X_\mu \ar[ld] \ar[rd] & \\ X_\mu & & Z'_\mu,}\] where $X_\mu=G({\bar k[[u]]}) s_\mu G({\bar k[[u]]})/\ci$ is a closed subscheme of $\Fl$, $\tilde X_\mu=G({\bar k[[u]]}) s_\mu G({\bar k[[u]]})/\ck_1$, and $Z'_\mu$ is the codimension-one $G_{\bar k} \times G_{\bar k}$-orbit in $X$ corresponding to $\mu$. The two maps in the diagram are smooth morphisms with isomorphic fibers. 

The geometric special fiber $M^{\loc} \otimes_{\co_E} \bar k$ is a closed reduced subscheme of $X_\mu$. And there exists a closed reduced subscheme $A'$ of $Z'_\mu$ such that the inverse image of $M^{\loc} \otimes_{\co_E} \bar k$ in $Z_\mu$ equals the inverse image of $A'$ in $Z_\mu$. Hence $M^{\loc} \otimes_{\co_E} \bar k$ is Cohen-Macaulay if and only if $A'$ is Cohen-Macaulay. 

By the explicit description of $\bar B$ in $X$ obtained by Brion \cite{B} and Springer \cite{Sp} and the description of $\mu$-admissible sets obtained in a joint work with Lam \cite{HL}, we'll show that $A'=\bar B \cap Z'_\mu$ is an open subscheme of the boundary of $B$ in $X$. By a result of Brion and Polo \cite{BP}, the boundary $\partial \bar B$ is Cohen-Macaulay. Hence $A'$ is Cohen-Macaulay. We finally obtain the Cohen-Macaulayness of the special fiber of the local model. 

In \cite{Z}, Zhu introduced global Schubert varieties, which are the (generalized) equal characteristic counterparts of the local models. It is also worth mentioning that by a similar argument, in the unramified case, global Schubert varieties associated to minuscule coweights are normal and Cohen-Macaulay. It would be interesting to see if it is still the case for arbitrary coweights. 

There is a different connection between local models and complete symmetric varieties, by Faltings in \cite{FF} and by Pappas (unpublished notes, see also \cite[Chapter 8]{PRS}). This approach doesn't use loop groups and works when the level subgroup is close to maximal parahoric. It would be interesting to compare the construction in this paper with their approach. 

It is also worth mentioning that in many cases, local charts around points of local models can be described by relatively simple matrix equations. Thus our results on local models imply structure results on some matrix equations. For instance, Theorem \ref{main}, together with an observation of G\"ortz \cite[Page 690]{G1}, implies that the equations $$A_1 A_2 A_3 \cdots A_n=A_2 A_3 \cdots A_n A_1=A_3 \cdots A_n A_1 A_2=\cdots=0,$$ where $A_1, \cdots, A_n$ are $r \times r$-matrices, define a Cohen-Macaulay singularity. 

\section{Local models}

\subsection{} In this section, we recall the definition and some results on $M^{\loc}$ in \cite{PZ}. 

Let $\co$ be a discrete valuation ring with fraction field $F$ and perfect residue field $k$ of characteristic
$p>0$. Fix a uniformizer $\varpi$ of $\co$. Let $G$ be a connected reductive group over $F$, split over a tamely ramified extension $\tilde F$ of $F$. Let $K$ be the parahoric subgroup of $G(F)$ associated to a vertex $x$ in the Bruhat-Tits building of $G(F)$. Let $\cg$ be the group scheme associated to $x$ in the sense of \cite[Theorem 3.2]{PZ}. It is a smooth affine group scheme over $\mathbb A^1_\co=\Spec(\co[u])$, the base change $\cg_{|\co[u^{\pm 1}]} \otimes_{\co[u^{\pm 1}]} F, u \mapsto \varpi$ is isomorphic to $G$ and the base change $\cg \otimes_{\co[u]} \co, u \mapsto \varpi$ is the parahoric group scheme of $G$ associated to $x$. 

For any $\co$-algebra $R$, we denote by $R[[u]]$ the ring of formal power series and $R((u))$ the ring of formal Laurent power series. We set $L\cg(R)=\cg(R((u)))$ and $L^+\cg(R)=\cg(R[[u]])$. Then $L\cg$ is represented by an ind-affine scheme over $\co$ and $L^+\cg$ is represented by an affine scheme over $\co$. Let $Gr_{\cg}=L\cg/L^+\cg$ be the fpqc quotient, which is represented by an ind-proper ind-scheme over $\co$. This is the {\it local affine Grassmannian}. See \cite[Proposition 6.3]{PZ}. 

\subsection{} 
Let $\{\mu\}$ be a geometric conjugacy class of one parameter subgroups of $G$. Let $E$ be the field of definition of $\{\mu\}$ and $E'=E F_{un}$, where $F_{un}$ is the maximal unramified extension of $F$ in $\tilde F$. As explained in \cite[7.a]{PZ}, there exists a representative of $\{\mu\}$ defined over $E'$ and this representative gives rise to an element $s_\mu$ in $L\cg(E')$. Moreover, the $L^+\cg$-orbit $(L^+\cg)_{E'} \cdot [s_\mu]$ in the affine Grassmannian $L\cg/L^+\cg \times_F E'$ is actually defined over $E$. In other words, there is an $E$-subvariety $X_\mu$ of $L\cg/L^+\cg \times_F E$ such that $X_\mu \times_E E'=(L^+\cg)_{E'} \cdot [s_\mu]$. 

The generalized local model $M_{\cg, \mu}$ (in the sense of Pappas and Zhu) is the reduced scheme over $\Spec(\co_E)$ which underlies the Zariski closure of $X_\mu$ in the ind-scheme $Gr_{\cg, \co_E}$. 

\subsection{} Let $\bar k$ be an algebraic closure of $k$ and $F'=\bar k((u))$. Let $G'=\cg \times_{\Spec(\co[u])} \Spec(F')$ be the base changing of $\cg$ to $F'$. Then $G'$ splits over a tamely ramified extension $\tilde F'$ of $F'$. Let $T'$ be the centralizer of a maximal split torus of $G'$. Let $I={\rm Gal}(\tilde F'/F')$ and $X_*(T')_I$ be the coinvariants of the coweight lattice $X_*(T')$. Let $\tW$ be the Iwahori-Weyl group of $G'$ and $W_0=N'(F')/T'(F')$ be the relative Weyl group of $G'$ over $F'$. There is a short exact sequence $$1 \to X_*(T')_I \to \tilde W \to W_0 \to 1.$$ 

For $\l \in X_*(T)_I$, we denote by $t_\l$ the corresponding translation element in $\tW$.\footnote{Here we adopt the sign convention in \cite{PZ}. In fact $t_\l$ equals $t^{\l}$ in \cite{HL}.}

To the geometric conjugacy class $\{\mu\}$ of one parameter subgroups of $G$, we associate a dominant coweight and denote it by $\mu$. Let $\L$ be the $W_0$-orbit in $X_*(T')_I$ that contains the image of $\mu$. Define the $\mu$-admissible set by $$\Adm(\mu)=\{w \in \tW; w \le t_\l \text{ for some } \l \in \L\}.$$ Here $\le$ is the Bruhat order on $\tilde W$ (see \cite[8.d.1]{PZ}). 

\subsection{}\label{adm} Recall that $x$ is a point in the Bruhat-Tits building of $G(F)$. Let $\cp'_x \subset G'$ be the corresponding parahoric group scheme over $\bar k[[u]]$. We choose a rational Borel $B'$ of $G'$ containing $T'$ in such a way that $\cp'_x$ is a standard parahoric group. 


Set $$\ca^{\cp'_x}(\mu)=\cup_{w \in \Adm(\mu)} L^+\cp'_x w L^+ \cp'_x/L^+\cp'_x.$$ It is a closed subscheme of the affine Grassmannian $Gr_{\cp'_x}$. 

The following result on the special fiber $M_{\cg, \mu} \otimes_{\co_E} k_E$ of the local model is obtained by Pappas and Zhu in \cite[Theorem 8.4 \& 8.5]{PZ}. 

\begin{thm}\label{pz}
Suppose that $p$ does not divide the order of the fundamental group of the derived group $\pi_1(G_{der})$. Then the special fiber $M_{\cg, \mu} \otimes_{\co_E} k_E$  is reduced and each geometric irreducible component is normal and Cohen-Macaulay. Moreover, the geometric special fiber $M_{\cg, \mu} \otimes_{\co_E} \bar k=\ca^{\cp'_x}(\mu)$ as closed subschemes of $Gr_{\cp'_x}$.
\end{thm}

Here the condition on $p$ is necessary to ensure that the corresponding loop group and affine Grassmannian variety are reduced. See \cite[Remark 6.4]{PR3}. 

In the case of a unitary or symplectic group that splits over an unramified extension, the local model $M_{\cg, \mu}$ coincides with the ``naive local model'' of Rapoport-Zink \cite{RZ} and the above properties of the special fiber were known earlier, by the work of G\"ortz \cite{G1} and \cite{G2}. 

\

The main purpose of this paper is to show that if $G$ splits over an unramified extension, then the special fiber of the local model, as a whole, is Cohen-Macaulay. By Serre's normality criterion and the behavior of depth under flat morphisms, this implies that\footnote{I am informed that Pappas and Zhu have recently proved the normality of local model (also in the non-split case) without using the Cohen-Macaulayness of the special fiber. }

\begin{thm}\label{main}
Suppose that $G$ splits over an unramified extension of $F$ and $p$ does not divide the order of the fundamental group of the derived group $\pi_1(G_{der})$. If $K$ is contained in a hyperspecial maximal compact subgroup of $G(F)$, then $M_{\cg, \mu}$ is normal and Cohen-Macaulay. 
\end{thm}

\section{Loop group and wonderful compactification}

\subsection{} In this section, we assume that $G$ is split over $F$. Hence $G'$ is also split over $F'$, i.e., $G'=LH$ is the loop group for some connected reductive algebraic group $H$ over $\bar k$. Let $T$ be a maximal torus of $H$ and $B \supset T$ be a Borel subgroup of $H$ such that $T'=T(F')$ and $B'=B(F')$. The pair $(B, T)$ determines the set of simple roots, which we denote by $S$. For any $J \subset S$, let $P_J \supset B$ be the standard parabolic subgroup of type $J$ and $P^-_J$ the opposite parabolic subgroup. Then $L_J=P_J \cap P^-_J$ is a standard Levi subgroup of $H$. For any parabolic subgroup $P$ of $H$, we denote by $U_P$ its unipotent radical. 

\subsection{}\label{zj} Now we recall the variety $Z_J$ introduced by Lusztig in \cite{Lu}. 

Let $J \subset I$. We define the action of $P^-_J \times P_J$ on $H \times H \times L_J$ by $(q, p) \cdot (h, h', z)=(h q \i, h' p \i, \pi_{P^-_J}(q) z \pi_{P_J}(p) \i)$. Here $\pi_{P^-_J}: P^-_J \to P^-_J/U_{P^-_J} \cong L_J$ and $\pi_{P_J}: P^-_J \to P_J/U_{P_J} \cong L_J$ are the projection maps. This is a free action. We denote by $Z_J=(H \times H) \times_{P^-_J \times P_J} L_J$ its quotient variety.

For any $h, h' \in H$ and $l \in L_J$, we denote by $[h, h', l]$ the image of $(h, h', l)$ in $Z_J$. The $H \times H$-action on $Z_J$ is defined by $(h, h') \cdot [a, b, c]=[h a, h' b, c]$. We write $h_J=[1, 1, 1]$. This is the base point of $Z_J$. 

\subsection{} Let $\pi: L^+H \to H$ be the reduction modulo $\varpi$ map and $\ck_1$ be the kernel of $\pi$. Then $\ck_1$ is a normal subgroup of $L^+H$. We define an action of $L^+H \times L^+H$ on $LH/\ck_1$ by $(h, h') \cdot z \ck_1=h z (h') \i \ck_1$. 

For any dominant coweight $\l \in X_*(T)$, set $$I(\l)=\{i \in S; \<\l, \a_i\>=0\}$$ and $$\tilde X_\l=L^+H s_\l L^+H/\ck_1,$$ where $s_\l$ is defined in $\S$1.2. 

Then $\tilde X_\l$ is a single $L^+H \times L^+H$-orbit and is a locally closed subscheme of $LH/\ck_1$. Moreover $LH/\ck_1=\sqcup_{\g} \tilde X_{\g}$, where $\g$ runs over all the dominant coweigths. 

The following result provides a relation between $\tilde X_\l$ and $Z_J$. 

\begin{prop}\label{11}
Let $\l$ be a dominant coweight and $J=I(\l)$. Then the map $L^+H \times L^+H \to Z_J$, $(h, h') \mapsto [\pi(h), \pi(h'), 1]=(\pi(h), \pi(h')) \cdot h_J$ induces a surjective $L^+H \times L^+H$-equivariant smooth morphism $$s: \tilde X_{\l} \to Z_J$$ and each fiber is isomorphic to an affine space over $\bar k$ of dimension $\<\l, 2 \rho\>-\ell(w_S)$. Here the action of $L^+H \times L^+H$ on $Z_J$ factors through the action of $H \times H$ on $Z_J$ defined in $\S$\ref{zj}, $\rho$ is the sum of all fundamental weights of $H$ and $w_S$ is the maximal element of $W_0$. 

\end{prop}

\begin{rmk}
An analogous result in mixed characteristic case is proved in a joint work with Wedhorn \cite{HW}. 
\end{rmk}

Proof. We first prove that $s$ is well-defined. We regard $H$ as a subgroup of $L^+H$. Then $L^+H=H \ck_1=\ck_1 H$. Since the map $L^+H \times L^+H \to Z_J$ is $H \times H$-equivariant, it suffices to show that

(a) For $h, h' \in H$ with $h s_\l (h') \i \in \ck_1 s_\l \ck_1$, $(h, h') \cdot h_J=h_J$. 

By assumption, $\emptyset \neq \ck_1 h \cap s_\l \ck_1 h' s_{-\l} \subset L^+H \cap s_\l L^+H s_{-\l}$. 

By \cite[Theorem 2.8.7]{Ca}, 
\begin{align*} L^+H \cap s_\l L^+H s_{-\l} &=(\ck_1 \cap s_\l \ck_1 s_{-\l}) (\ck_1 \cap s_\l H s_{-\l}) (H \cap s_\l \ck_1 s_{-\l}) (H \cap s_\l H s_{-\l}).
\end{align*}

We have that $\ck_1 \cap s_\l H s_{-\l}=s_{\l} U_{P_J} s_{-\l}$, $H \cap s_\l \ck_1 s_{-\l}=U_{P_J^-}$ and $H \cap s_\l H s_{-\l}=L_J$. Then there exists $z \in \ck_1 \cap s_\l \ck_1 s_{-\l}$, $l \in L_J$, $u \in U_{P_J^-}$ and $u' \in U_{P_J}$ such that $$z (s_{\l} u' s_{-\l}) u l \in \ck_1 h \cap s_\l \ck_1 h' s_{-\l}.$$

Notice that $s_{\l} u' s_{-\l} \in \ck_1$. Hence $z (s_{\l} u' s_{-\l}) u l \in \ck_1 u l$ and $h=u l$. Similarly, $s_{-\l} u s_\l \in \ck_1$ and $s_{-\l} \bigl(z (s_{\l} u' s_{-\l}) u l \bigr) s_{\l} \in \ck_1 u' l$. Hence $h'=u' l$. Therefore $$(h, h') \cdot h_J=(u l, u' l) \cdot h_J=(u, u') \cdot h_J=h_J.$$

(a) is proved. 

Since $L^+H \times L^+H$ acts transitively on $\tilde X_{\l}$ and on $Z_J$, and the map $s: \tilde X_{\l} \to Z_J$ is $L^+H \times L^+H$-equivariant, all fibers are isomorphic. Notice that $Z_J$ is reduced. Thus $s$ is generically flat and hence is flat by equivariance. 

Now we consider the fiber over $h_J$. It is $$\{\ck_1 U_{P^-_J} l s_\l l \i U_{P_J} \ck_1; l \in L_J\}/\ck_1=\ck_1 U_{P^-_J} s_\l U_{P_J} \ck_1/\ck_1.$$

Since $s_{-\l} U_{P^-_J} s_\l \subset \ck_1$ and $s_{\l} U_{P_J} s_{-\l} \subset \ck_1$, we have that 
\begin{align*}
\ck_1 U_{P^-_J} s_\l  U_{P_J} \ck_1 &=\ck_1 s_\l (s_{-\l} U_{P^-_J} s_\l)  U_{P_J} \ck_1 \subset \ck_1 s_\l \ck_1 U_{P_J} \ck_1 \\ &=\ck_1 s_\l U_{P_J} \ck_1=\ck_1 (s_\l U_{P_J} s_{-\l}) s_\l \ck_1 \\ & \subset \ck_1 s_\l \ck_1. 
\end{align*} 
It is obvious that $\ck_1 s_\l \ck_1 \subset \ck_1 U_{P^-_J} s_\l  U_{P_J} \ck_1$. Therefore $\ck_1 U_{P^-_J} s_\l  U_{P_J} \ck_1=\ck_1 s_\l \ck_1$ and 
$$\ck_1 U_{P^-_J} s_\l  U_{P_J} \ck_1/\ck_1=\ck_1 s_\l \ck_1/\ck_1 \cong \ck_1/(\ck_1 \cap s_{\l} \ck_1 s_{-\l}).$$
This is an affine space of dimension $\dim(\tilde X_\l)-\dim(Z_J)=\dim(X_\l)+\dim(B)-\dim(G)=\<\l, 2 \rho\>-\ell(w_S)$. 
\qed




\subsection{} Now we recall the definition and some elementary facts on the wonderful compactification. More details can be found in the survey article of Springer \cite{Sp1}. 

Let $\underline H$ be the adjoint group of $H$. The set of simple roots of $\underline H$ is again denoted by $S$. For any subgroup $H'$ of $H$, we denote by $\underline H'$ the image of $H'$ via the map $H \to \underline H$. 

Let $X$ be the wonderful compactification of $\underline H$ (\cite{CP}, \cite{Str}). Roughly speaking, one starts with a suitable finite-dimensional projective representation $\rho: \underline H \to PGL(V)$ of $H$, then $X$ is defined to be the closure in $\mathbb P(\End(V))$ of the image $\rho(\underline H)$. The closure is independent of the choice of $\rho$. 

It is known that $X$ is an irreducible, smooth projective $(\uh \times \uh)$-variety with finitely many $\uh \times \uh$-orbits indexed by the subsets of $S$. They are described as follows. 

Let $J \subset S$. Let $\uh_J$ be the adjoint group of $\underline L_J$ (and hence of $L_J$). We define an action of $\underline P^-_J \times \underline P_J$ on $\underline H \times \underline H \times \uh_J$ in the same way as in $\S$\ref{zj} and denote by $Z'_J=(\underline H \times \underline H) \times_{\underline P^-_J \times \underline P_J} \uh_J$ the quotient variety. The group $\underline H \times \underline H$ acts on $Z'_J$ in the same way as in $\S$\ref{zj}. We denote by $h'_J$ the image in $Z'_J$ of $(1, 1, 1) \in \underline H \times \underline H \times \uh_J$. This is the base point of $Z'_J$. It is known that $X=\sqcup_{J \subset S} Z'_J$ as the union of $\uh \times \uh$-orbits. 

For any locally closed subscheme $Z$ of $X$, we denote by $\bar Z$ the closure of $Z$ in $X$. The closure relation between $\uh \times \uh$-orbits on $X$ is described as follows. For any $J \subset S$, $$\bar Z'_J=\sqcup_{J' \subset J} Z'_{J'}.$$ In particular, $Z_S=\underline H$ is the open orbit in $X$ and for any maximal proper subset $J$ of $S$, $Z'_J$ is a codimension-one orbit of $X$ and hence is open in the boundary $\partial X=X \backslash \underline H$ of $\underline H$. The closed orbit $Z'_\emptyset$ is isomorphic to $\underline H/\ub \times \underline H/\ub$. 

\subsection{}\label{corr} Now we discuss the situation we may apply to the study of local models. By the choice of $B$ and $T$, $L^+\cp'_x$ is the inverse image of $P_Y$ under $\pi$ for some $Y \subset S$ and $L^+\cp'_x/\ck_1 \cong P_Y$. Therefore the projection map $f_1: LH/\ck_1 \to LH/L^+\cp'_x$ is a $P_Y$-torsor. Hence the restriction $f_1: \tilde X_{\mu} \to X_{\mu}$ is a $P_Y$-torsor. 

Now we have the following diagram \[\xymatrix{X_{\mu} & \tilde X_{\mu} \ar[l]_-{f_1} \ar[r]^-s & Z_{I(\mu)} \ar[r]^-{f_2} & Z'_{I(\mu)}}.\] Here $f_2: Z_{I(\mu)} \to Z'_{I(\mu)}$ is induced from the map $H \times H \times L_{I(\mu)} \to \underline H \times \underline H \times \uh_{I(\mu)}$ and hence is a smooth morphism with fibers isomorphic to the center of $L_{I(\mu)}$. 

Notice that $s, f_1, f_2$ are smooth morphisms. Hence if $A$ is a closed reduced subscheme of $X_{\mu}$ and $A'$ is a closed reduced subscheme of $Z'_{I(\mu)}$ such that $f_1 \i(A)=(f_2 \circ s)  \i(A')$, then $A$ is Cohen-Macaulay if and only if $A'$ is Cohen-Macaulay. 

\section{boundary of parabolic subgroup}

In this section, we study the boundary in $X$ of the parabolic subgroup $\underline P_Y$ of $\underline H$. As we'll see in the next section, this boundary is closely related to the geometric special fiber of the local model in the sense of $\S$\ref{corr}. 

\subsection{} We first recall some results on the $\underline B \times \underline B$-orbits of $X$ obtained by Brion \cite{B} and Springer \cite{Sp}. 

For any $J \subset S$, we denote by $W_J$ the subgroup of $W_0$ generated by the simple reflections in $J$ and $W^J_0$ the set of minimal length representatives in $W_0/W_J$. Let $w_J$ be the maximal element in $W_J$. 

For any $(x, y) \in W^J_0 \times W_0$, we set $$[J, x, y]=(\underline B x, \underline B y) \cdot h'_J \subset Z'_J.$$

By \cite[Lemma 1.3]{Sp}, $X=\sqcup_{J \subset S} \sqcup_{x \in W^J_0, y \in W_0} [J, x, y]$. 

The closure relations between $\ub \times \ub$-orbits of $X$ has been obtained in \cite[Proposition 2.4]{Sp}. The following simplified version is found in \cite[Proposition 6.3]{HT}. 

\begin{prop}\label{border}
Let $J, J' \subset S$, $x \in W^J_0$, $x' \in W^{J'}$ and $y, y' \in W_0$. Then $[J', x', y'] \subset \overline{[J, x, y]}$ if and only if $J' \subset J$ and there exists $u \in W_J$ such that $x u \le x'$, $y' \le y u$. 
\end{prop}

\subsection{}\label{special} Now we discuss some special cases that will be used in this paper. 

For $J \subset S$, we define a partial order $\preceq_J$ on $W^J_0 \times W_0$ as follows. Let $(x, y), (x', y') \in W^J_0 \times W_0$, we write $(x', y') \preceq_J (x, y)$ if there exists $u \in W_J$ such that $x u \le x'$, $y' \le y u$. Then $[J, x', y'] \subset \overline{[J, x, y]}$ if and only if $(x', y') \preceq_J (x, y)$. 

The following joint result with Lam \cite[Theorem 2.2]{HL} relates this partial order with the Bruhat order on the Iwahori-Weyl group. 

\begin{prop}\label{hl} Let $\mu$ be a minuscule coweight. Then 

(1) The map $$(W^{I(\mu)}_0 \times W_0, \preceq_{I(\mu)}) \to (W_0 t_\mu W_0, \le), \quad (x, y) \mapsto x t_\mu y \i$$ is a bijection between posets. Here $\le$ is the restriction to $W_0 t_\mu W_0$ of the Bruhat order on $\tW$. 

(2) Set $Q_\mu=\{(x, y) \in W^{I(\mu)}_0 \times W_0; y \le x\}$. Then the restriction of the map in (1) gives a bijection from $Q_\mu$ to the admissible set $\Adm(\mu)$. 
\end{prop}

\begin{rmk}
By definition, the maximal elements of $Q_\mu$ are $(x, x)$ for $x \in W^{I(\mu)}_0$. The bijection in (1) send these elements to the elements $t_\l$ for $\l$ in the $W_0$-orbit of $\mu$, which are just the maximal elements in $\Adm(\mu)$. 
\end{rmk}



\subsection{} As a special case of Proposition \ref{border}, the closure of $\up_Y$ in $X$ is described as follows \[\overline{\up_Y}=\overline{[S, 1, w_Y]}=\sqcup_{J \subset S}\sqcup_{x \in W^J_0, y \in W_0, \min(W_Y y) \le x} [J, x, y].\]

For our purpose, we need a different description of $\overline{\up_Y}$. 

\begin{cor}\label{upy}
For any $J \subset S$, $$\overline{\up_Y} \cap Z'_J=\cup_{w \in W^J_0} \overline{(\up_Y w, \up_Y w) \cdot h'_J} \cap Z'_J.$$
\end{cor}

Proof. By Proposition \ref{border} and $\S$\ref{special}, $$\overline{\ub} \cap Z'_J=\sqcup_{(x, y) \in W^J_0 \times W_0, y \le x} [J, x, y] \subset \cup_{w \in W^J_0} \overline{[J, w, w]} \cap Z'_J.$$ On the other hand, $\overline{\ub} \cap Z'_J$ is closed in $Z'_J$ and $[J, w, w] \subset \overline{\ub}$ for all $w \in W^J_0$. Therefore $\overline{\ub} \cap Z'_J=\cup_{w \in W^J_0} \overline{[J, w, w]} \cap Z'_J$.

As $\overline{\ub}$ is closed in $X$ and stable by $\ub \times \ub$, and $\ub$ is a Borel subgroup of $\up_Y$, $(\up_Y, \up_Y) \cdot \overline{\ub}$ is closed in $X$ and thus equals $\overline{\up_Y}$. Similarly, for any $w \in W^J_0$, $(\up_Y, \up_Y) \cdot (\overline{[J, w, w]} \cap Z'_J)$ is closed in $Z'_J$ and equals $\overline{(\up_Y w, \up_Y w) \cdot h'_J} \cap Z'_J$. Hence  
\begin{align*} \overline{\up_Y} \cap Z'_J &=(\up_Y, \up_Y) \cdot \overline{\ub} \cap Z'_J=(\up_Y, \up_Y) \cdot (\overline{\ub} \cap Z'_J) \\ & =\cup_{w \in W^J_0} (\up_Y, \up_Y) \cdot (\overline{[J, w, w]} \cap Z'_J) \\&=\cup_{w \in W^J_0} \overline{(\up_Y w, \up_Y w) \cdot h'_J} \cap Z'_J.\end{align*}\qed


\subsection{} It is proved by Brion and Polo in \cite[Theorem 20]{BP} that $\overline{\up_Y}$ is Cohen-Macaulay. Notice that $\partial X$ intersects properly $\overline{\up_Y}$. Let $\partial \overline{\up_Y}=\overline{\up_Y}  \cap \partial X$ be the boundary of $\overline{\up_Y}$ in $X$. Since $\partial \overline{\up_Y}$ is a hypersurface of $\overline{\up_Y}$, with local equation being a nonzero divisor, we have that

\begin{prop}
$\partial \overline{\up_Y}$ is Cohen-Macaulay. 
\end{prop}

\section{Proof of the main theorem} 



As local models are compatible with unramified base change, it suffices to consider the case where $G$ is split over $F$. We keep this assumption in the rest of this section. 

\subsection{} The geometric special fiber $M_{\cg, \mu} \otimes_{\co_E} \bar k=\ca^{\cp'_x}(\mu)$ is a closed subscheme of $X_\mu$. Then $\tilde \ca^{\cp'_x}(\mu)=\cup_{w \in \Adm(\mu)} L^+\cp'_x w L^+ \cp'_x/\ck_1$ is a reduced closed subscheme of $\tilde X_\mu$ and is the inverse image of $\ca^{\cp'_x}(\mu)$ under the map $f_1$. 

Let $A^Y(\mu)$ be the reduced subscheme of $Z'_{I(\mu)}$, which equals 
$$\cup_{(x, y) \in Q_\mu} (\up_Y x, \up_Y y) \cdot h'_{I(\mu)}=\cup_{x \in W^{I(\mu)}_0} \overline{(\up_Y x, \up_Y x) \cdot h'_{I(\mu)}} \cap Z'_{I(\mu)}$$ as a set. Notice that $\cp'_x=P_Y \ck_1=\ck_1 P_Y$. By Proposition \ref{hl} (2), we have that $\tilde \ca^{\cp'_x}(\mu)=(f_2 \circ s) \i A^Y(\mu)$. 

Thus the reduced schemes $\ca^{\cp'_x}(\mu)$ and $A^Y(\mu)$ are related in the sense of $\S$\ref{corr}. 

\subsection{} We'll then prove that $A^Y(\mu)$ is the scheme-theoretic intersection of $\overline{\up_Y}$ with $Z'_{I(\mu)}$. 

We first show that $A^Y(\mu)=\overline{\up_Y} \cap Z'_{I(\mu)}$ set-theoretically. 

By definition, $\ca^{\cp'_x}(\mu)$ is the union of the closures of $L^+\cp'_x s_{w \mu} L^+\cp'_x/L^+\cp'_x$ in $LH/L^+\cp'_x$, where $w$ runs over elements in $W^{I(\mu)}_0$. By $\S$\ref{corr}, $A^Y(\mu)$ is the union of the closures of $(\up_Y w, \up_Y w) \cdot h'_{I(\mu)}$ in $Z'_{I(\mu)}$, where $w$ again runs over elements in $W^{I(\mu)}_0$. Hence by Corollary \ref{upy}, $A^Y(\mu)=\overline{\up_Y} \cap Z'_{I(\mu)}$ as sets. 

It remains to show that $\overline{\up_Y} \cap Z'_{I(\mu)}$ is reduced. We recall a result in \cite[Proposition 6.2]{HT}, which strengthened \cite[Theorem 2]{BP}. 

\begin{prop}
There exists a Frobenius splitting on $X$ that compatibly splits all the $\ub \times \ub$-orbit closures. 
\end{prop}

In particular, there exists a Frobenius splitting on $X$ that compatibly splits $\overline{\up_Y}$ and $\overline{Z'_{I(\mu)}}$. Therefore the scheme-theoretic intersection $\overline{\up_Y} \cap \overline{Z'_{I(\mu)}}$ is a split scheme and hence is reduced. So the scheme-theoretic intersection $\overline{\up_Y} \cap Z'_{I(\mu)}$ is also reduced. 

\subsection{} Now we prove Theorem \ref{main} for the split case. 

Since $\mu$ is minuscule, $I(\mu)$ is a maximal proper subset of $S$. Hence $Z'_{I(\mu)}$ is an open subscheme of $\partial X$ and $A^Y(\mu)=\overline{\up_Y} \cap Z'_{I(\mu)}$ is an open subscheme of $\partial \overline{\up_Y}$. Since $\partial \overline{\up_Y}$ is Cohen-Macaulay, $A^Y(\mu)$ is also Cohen-Macaulay. 

By $\S$\ref{corr}, the geometric special fiber $\ca^{\cp'_x}(\mu)$ is Cohen-Macaulay and so is the special fiber $M_{\cg, \mu} \otimes_{\co_E} k_E$. Since $M_{\cg, \mu}$ is flat over $\text{Spec} (\co_E)$, it is Cohen-Macaulay. 

By Theorem \ref{pz}, $M_{\cg, \mu} \otimes_{\co_E} k_E$ is generically smooth. By Serre's criterion, $M_{\cg, \mu}$ is also normal. This finishes the proof. 

\section*{Acknowledgement} We thank U. G\"ortz, T. Wedhorn and X. Zhu for helpful discussions on local models and J. Starr for answering my questions on algebraic geometry. We thank G. Pappas for pointing out a different connection between local models and wonderful compactifications \cite{FF}, \cite[Chapter 8]{PRS}. We also thank the referees for many useful suggestions.

\end{document}